\documentclass[12pt,reqno]{amsart}

\usepackage{amssymb}

\newtheorem{theorem}{Theorem}[section]
\newtheorem{proposition}[theorem]{Proposition}
\newtheorem{lemma}[theorem]{Lemma}
\newtheorem{corollary}[theorem]{Corollary}

\numberwithin{equation}{section}

\begin{document}

\baselineskip=16pt

\title{On connections on principal bundles}

\author[I. Biswas]{Indranil Biswas}

\address{School of Mathematics, Tata Institute of Fundamental
Research, Homi Bhabha Road, Mumbai 400005, India}

\email{indranil@math.tifr.res.in}

\date{}



\begin{abstract}
A new construction of a universal connection was given in \cite{BHS}. The main
aim here is to explain this construction. A theorem of Atiyah and Weil says that
a holomorphic vector bundle $E$ over a compact Riemann surface admits a holomorphic connection
if and only if the degree of every direct summand of $E$ is degree. In \cite{AB}, this
criterion was generalized to principal bundles on compact Riemann surfaces. This
criterion for principal bundles is also explained.
\end{abstract}

\maketitle

\section{Introduction}\label{sec0}

A connection $\nabla^0$ on a $C^\infty$ principal $G$--bundle ${\mathcal E}_G\, \longrightarrow\,
{\mathcal X}$ is called \textit{universal} if given any $C^\infty$ principal $G$--bundle $E_G$
on a finite dimensional $C^\infty$ manifold $M$, and any connection $\nabla$ on $E_G$, there
is a $C^\infty$ map
$$
\xi\,:\, M\, \longrightarrow\, \mathcal X
$$
such that
\begin{itemize}
\item the pulled back principal $G$--bundle $\xi^*{\mathcal E}_G$ is isomorphic to $E_G$, and

\item the isomorphism between $\xi^*{\mathcal E}_G$ and $E_G$ can be so chosen that
it takes the pulled back connection $\xi^*\nabla^0$ on $\xi^*{\mathcal E}_G$ to the connection
$\nabla$ on $E_G$.
\end{itemize}
In \cite{NR} and \cite{Sc} universal connections were constructed. In \cite{BHS} a very simple,
in fact quite tautological, universal connection was constructed.

\section{Atiyah bundle}\label{sec1}

All manifolds considered here will be $C^\infty$, 
second countable and Hausdorff. Later we will impose further conditions
such as complex structure.

Let $G$ be a finite dimensional Lie group. Take a connected $C^\infty$ manifold $M$.
A principal $G$--bundle over $M$ is a triple of the form
\begin{equation}\label{tri}
(E_G,\,  p,\, \psi)\, ,
\end{equation}
where
\begin{enumerate}
\item{} $E_G$ is a $C^\infty$ manifold manifold,

\item{}
\begin{equation}\label{e1}
p\, :\, E_G\, \longrightarrow\, M
\end{equation}
is a $C^\infty$ surjective
submersion, and

\item{}
\begin{equation}\label{e2}
\psi\, :\, E_G\times G\, \longrightarrow\, E_G
\end{equation}
is a $C^\infty$ map defining a right action of $G$ on $E_G$, such that the
following two conditions hold:
\end{enumerate}
\begin{itemize}
\item{} the two maps $p\circ\psi$ and $p\circ p_1$
from $E_G\times G$ to $M$ coincide, where $p_1$ is the natural
projection of $E_G\times G$ to $E_G$, and
\item{} the map to the fiber product
$$
\text{Id}_{E_G}\times\psi\, :\, E_G\times G\, \longrightarrow\,
E_G\times_M E_G
$$
is a diffeomorphism; note that the first condition $p\circ\psi\, =\, p\circ p_1$
implies that the image of $\text{Id}_{E_G}\times\psi$ is contained
in the submanifold $E_G\times_M E_G\, \subset\, E_G\times E_G$ consisting of all point $(z_1,\, z_2)
\,\in\, E_G\times E_G$ such that $p(z_1)\, =\, p(z_2)$.
\end{itemize}
Therefore, the first condition implies that $G$ acts on $E_G$ along the fibers of $p$,
while the second condition implies that the action of $G$ on each fiber of $p$
is both free and transitive.

Take a $C^\infty$ principal $G$--bundle $(E_G,\, p,\, \psi)$ over $M$.
The tangent bundle of the manifold $E_G$ will be denoted by $TE_G$.
Take a point $x\,\in\, M$. Let
$$
(TE_G)^x\,:=\, (TE_G)\vert_{p^{-1}(x)}\,\longrightarrow\,p^{-1}(x)
$$
be the restriction of the vector bundle
$TE_G$ to the fiber $p^{-1}(x)$ of $p$ over the point $x$. As noted above,
the action $\psi$
of $G$ on $E_G$ preserves $p^{-1}(x)$, and the resulting action
of $G$ on $p^{-1}(x)$ is free and transitive. Therefore, the action of $G$
on $TE_G$ given by $\psi$ restricts to an action of $G$ on
$(TE_G)^x$. Let $\text{At}(E_G)_x$ be the space of all $G$--invariant
sections of $(TE_G)^x$. Since the action of $G$ on the
fiber $p^{-1}(x)$ is transitive,
it follows that any $G$--invariant section of $(TE_G)^x$ is automatically smooth.
More precisely, any $G$--invariant sections of $(TE_G)^x$ is uniquely determined
by its evaluation of some fixed point of $p^{-1}(x)$. Therefore, $\text{At}(E_G)_x$
is a real vector space whose dimension coincides with the dimension of $E_G$.

There is a natural vector bundle over $M$, which was
introduced in \cite{At}, whose fiber over any
$x\,\in\, M$ is $\text{At}(E_G)_x$. This vector bundle is known
as the \textit{Atiyah bundle}, and it is denoted by
$\text{At}(E_G)$. We now recall the construction of $\text{At}(E_G)$.

As before, consider the action of $G$ to $TE_G$ given by the action $\psi$ of $G$ on $E_G$.
Since the action of $G$ is free and transitive on each fiber of $p$, it follows
that this action of $G$ on $TE_G$ free and proper. Therefore, we have a quotient manifold
\begin{equation}\label{e6}
\text{At}(E_G)\, :=\, (TE_G)/G
\end{equation}
for this action of $G$ on $TE_G$. Since the natural projection $TE_G\, \longrightarrow\, E_G$
is $G$--equivariant, it produces a projection
\begin{equation}\label{ap}
\text{At}(E_G)\, :=\, (TE_G)/G\, \longrightarrow\, E_G/G \,=\, M\, .
\end{equation}
This projection in \eqref{ap} is clearly surjective. Furthermore, it is a submersion because the
projection $TE_G\, \longrightarrow\, E_G$ is so. It is now straight--forward to check that the
projection in \eqref{ap} makes $\text{At}(E_G)$ a
$C^\infty$ vector bundle over $M$. Its rank coincides
with the rank of the tangent bundle $TE_G$, so its rank is $\dim G +\dim M$. From \eqref{e6}
it follows immediately that we have a natural diffeomorphism
\begin{equation}\label{i.}
\mu\, :\, p^*\text{At}(E_G)\, \longrightarrow\, TE_G\, .
\end{equation}
It is straight--forward to check that $\mu$ is a $C^\infty$ isomorphism of vector bundles
over $E_G$.

Let
\begin{equation}\label{dp}
dp \, :\, TE_G\, \longrightarrow\, p^*TM
\end{equation}
be the differential of the projection $p$ in
\eqref{e1}. Consider the surjective $C^\infty$ homomorphism
of vector bundles
\begin{equation}\label{dx}
dp\circ\mu\, :\, p^*\text{At}(E_G)\, \longrightarrow\, p^*TM\, ,
\end{equation}
where $\mu$ is constructed in \eqref{i.}.
Since $p^*\text{At}(E_G)$ and $p^*TM$ are
pulled back to $E_G$ from $M\,=\, E_G/G$,
they are naturally equipped with an action of $G$.
The homomorphism $dp\circ\mu$ in \eqref{dx} is clearly
$G$--equivariant. Therefore, it descends to a surjective
$C^\infty$ homomorphism of vector bundles
\begin{equation}\label{s.h.}
\eta\, :\, \text{At}(E_G)\, \longrightarrow\, TM\, .
\end{equation}

The kernel of the differential $dp$ in \eqref{dp} is clearly
preserved by the action of $G$ on $TE_G$. The quotient
$\text{kernel}(dp)/G$ will be denoted by $\text{ad}(E_G)$.
It is a $C^\infty$ vector bundle on $M$ whose rank is $\dim G$.
The inclusion of $\text{kernel}(dp)$ in $TE_G$ produces a
fiberwise injective $C^\infty$ homomorphism of vector bundles
$$
\iota_0\, :\, \text{ad}(E_G)\,\longrightarrow\, \text{At}(E_G)\, .
$$
The kernel of the homomorphism $\eta$ in \eqref{s.h.} coincides
with the image of $\iota_0$. Therefore, we have a short exact
sequence of $C^\infty$ vector bundles over $M$
\begin{equation}\label{e7}
0\,\longrightarrow\,{\rm ad}(E_G)\,
\stackrel{\iota_0}{\longrightarrow}
\, {\rm At}(E_G)\,\stackrel{\eta}{\longrightarrow}\, TM 
\,\longrightarrow\,0\, ,
\end{equation}
which is known as the \textit{Atiyah exact sequence} for $E_G$. Using the
Lie bracket operation of vector fields on $E_G$, the fibers of $\text{ad}(E_G)$
are Lie algebras; this will be elaborated below.

The Lie algebra of $G$ will be denoted by $\mathfrak g$.
Consider the action of $G$ on itself defined by
$\text{Ad}(g)(h)\,=\, g^{-1}hg$. This action
defines an action of $G$ on $\mathfrak g$,
which is known as the \textit{adjoint action}; this adjoint action 
of $G$ on $\mathfrak g$ will also be denoted by $\text{Ad}$.
Consider the quotient of $E_G\times
{\mathfrak g}$ where two points $(z, \,
v),\, (z',\, v')\, \in\, E_G\times{\mathfrak g}$ are identified
if there is some $g_0\, \in\, G$ such
that $z'\,=\,zg_0$ and $v'\,=\, \text{Ad}(g^{-1}_0)(v)$.
This quotient space coincides with the total space of the
adjoint vector bundle $\text{ad}(E_G)$ in \eqref{e7}. Note that
the projection
\begin{equation}\label{ad.d}
\text{ad}(E_G)\, \longrightarrow\, M\, .
\end{equation}
sends the
equivalence class of any $(z,\, v)\, \in\, E_G\times{\mathfrak g}$
to $p(z)$ (it is clearly independent of the choice of the element in
the equivalence class). The fibers of $\text{ad}(E_G)$ are
identified with $\mathfrak g$ up to conjugation. Since the adjoint
action of $G$ on $\mathfrak g$ preserves its Lie algebra structure,
the fibers of $\text{ad}(E_G)$ are in fact Lie algebras isomorphic to
$\mathfrak g$. This Lie algebra structure of a fiber of $\text{ad}(E_G)$ coincides with the
one constructed earlier using the Lie bracket operation of vector fields. The pulled back
vector bundle $p^*\text{ad}(E_G)$ on $E_G$ is identified with the trivial vector bundle
$E_G\times\mathfrak g$ with fiber $\mathfrak g$. This identification sends any vector
$(z,\, v)\, \in\, (p^*\text{ad}(E_G))_z$ in the fiber over $z$ of the pulled back bundle 
to the element $(z,\, v)$ of the trivial vector bundle $E_G\times\mathfrak g$.

A \textit{connection} on $E_G$ is a $C^\infty$ splitting of
the Atiyah exact sequence for $E_G$ \cite{At}. In other words,a connection on
$E_G$ is a $C^\infty$ homomorphism of vector bundles
\begin{equation}\label{D}
D\, :\, TM\,\longrightarrow\,{\rm At}(E_G)
\end{equation}
such that $\eta\circ D\, =\,\text{Id}_{TM}$, where $\eta$
is the projection in \eqref{s.h.}.

Let
\begin{equation}\label{d2}
D\, :\, TM\,\longrightarrow\,{\rm At}(E_G)
\end{equation}
be a homomorphism defining a connection on $E_G$.
Consider the composition homomorphism
$$
p^*TM\,\stackrel{p^*D}{\longrightarrow}\,p^*{\rm At}(E_G)
\,\stackrel{\mu}{\longrightarrow}\, TE_G\, ,
$$
where $\mu$ is the isomorphism in \eqref{i.}. Its image
\begin{equation}\label{e8}
{\mathcal H}(D)\, :=\, (\mu\circ p^*D) (p^*TM)\, \subset\, TE_G
\end{equation}
is known as the \textit{horizontal subbundle}
of $TE_G$ for the connection $D$. Since $\mu$
is an isomorphism, and the splitting homomorphism
$D$ in \eqref{d2} is uniquely determined by its image
$D(TM)\, \subset\, {\rm At}(E_G)$, it follows immediately that
the horizontal subbundle
${\mathcal H}(D)$ determines the connection $D$ uniquely.

The composition
$$
\text{kernel}(dp)\, \hookrightarrow\, TE_G
\,\longrightarrow\, TE_G/{\mathcal H}(D)
$$
is an isomorphism. Hence we have
$$
TE_G\, =\, {\mathcal H}(D) \oplus (E_G\times {\mathfrak g})\, ;
$$
it was noted earlier that $p^*\text{ad}(E_G)$ is identified with
the trivial vector bundle $E_G\times {\mathfrak g}$.
The projection of $TE_G$ to the second factor of the above
direct sum decomposition defines a $\mathfrak g$--valued
smooth one--form on $E_G$. The connection $D$ is clearly
determined uniquely by this $\mathfrak g$--valued one--form
on $E_G$.

See \cite[p.~370, Lemma~2.2]{BHS} for a proof of the following lemma:

\begin{lemma}\label{lem1}
Any principal $G$--bundle $E_G\, \longrightarrow\, M$ admits
a connection.

The space of all connections on a principal $G$--bundle $E_G$
is an affine space for
the vector space $C^\infty(M;\, {\rm Hom}(TM,\, {\rm ad}(E_G)))$.
\end{lemma}

\section{A universal connection}\label{sec4}

\subsection{A tautological connection}\label{sec3}

As before,
let $p\, :\, E_G\, \longrightarrow\, M$ be a $C^\infty$ principal $G$--bundle.
Consider the Atiyah exact sequence in \eqref{e7}.
Tensoring it with the cotangent bundle
$T^*M\, =\, (TM)^*$ we get the following
short exact sequence of vector bundles on $M$
\begin{equation}\label{e10}
0\,\longrightarrow\, {\rm ad}(E_G)\otimes T^*M\,\longrightarrow
\,{\rm At}(E_G)\otimes T^*M\, \stackrel{\eta\otimes 
\text{Id}_{T^*M}}{\longrightarrow}\, TM\otimes T^*M
\,=:\, {\rm End}(TM)\,\longrightarrow\, 0\, .
\end{equation}
Let $\text{Id}_{TM}$ denote the identity automorphism of
$TM$. It defines a $C^\infty$ section of the endomorphism bundle
${\rm End}(TM)$. Let
\begin{equation}\label{e11}
\delta\, :\, {\mathcal C}(E_G)\, :=\, (\eta\otimes 
\text{Id}_{T^*M})^{-1}(\text{Id}_{TM})\, \subset\, {\rm 
At}(E_G)\otimes T^*M \,\longrightarrow\, M
\end{equation}
be the fiber bundle over $M$, where $\eta\otimes
\text{Id}_{T^*M}$ is the surjective homomorphism in \eqref{e10}.

We recall that a connection on $E_G$ is a $C^\infty$ splitting of the
Atiyah exact sequence.

See \cite[p.~371, Lemma~3.1]{BHS} for a proof of the following:

\begin{lemma}\label{lem2}
The space of all connections on $E_G$ is in bijective
correspondence with the space of all smooth sections of the
fiber bundle
$$
\delta\, :\, {\mathcal C}(E_G)\, \longrightarrow\, M
$$
constructed in \eqref{e11}.
\end{lemma}

Combining Lemma \ref{lem1} with Lemma \ref{lem2}, the following is obtained.

\begin{corollary}\label{cor-2}
The fiber bundle $\delta$ in \eqref{e11} is an affine bundle
over $M$ for the vector bundle ${\rm Hom}(TM,\, {\rm ad}(E_G))$.
In particular, if we fix a connection on $E_G$ (which exists
by Lemma \ref{lem1}), then the fiber bundle in \eqref{e11}
gets identified with the total space of the vector bundle
${\rm Hom}(TM,\, {\rm ad}(E_G))$.
\end{corollary}

See \cite[p.~372, Proposition~3.3]{BHS} for a proof of the following:

\begin{proposition}\label{prop2}
There is a tautological connection on the principal
$G$--bundle $\delta^*E_G$ over ${\mathcal C}(E_G)$.
\end{proposition}

The key observations in the construction of the tautological connection in
Proposition \ref{prop2} are the following:

There is a tautological homomorphism
$$
\beta\, :\, \delta^*\text{At}(E_G)\, \longrightarrow\,
\delta^*\text{ad}(E_G)\,=\, \text{ad}(\delta^*E_G)\, .
$$
On the other hand, there is a tautological projection
$$
\beta'\, :\, \text{At}(\delta^* E_G)\, \longrightarrow\,
\delta^*\text{At}(E_G)$$
such that the diagram
$$
\begin{matrix}
\text{At}(\delta^* E_G)& \stackrel{\beta'}{\longrightarrow} &
\delta^*\text{At}(E_G)\\
\Big\downarrow && ~\,~\,~\,~\,\Big\downarrow \delta^*\eta\\
T{\mathcal C}(E_G)& \stackrel{d\delta}{\longrightarrow} &
\delta^*TM
\end{matrix}
$$
where the projection $\text{At}(\delta^* E_G)\, \longrightarrow\,
T{\mathcal C}(E_G)$ is constructed as in \eqref{s.h.} for
the principal $G$--bundle $\delta^* E_G$. Finally, the composition
$$
\beta\circ\beta'\, :\, \text{At}(\delta^* E_G)\, \longrightarrow\,
\text{ad}(\delta^* E_G)
$$
gives a splitting of the Atiyah exact sequence for $\delta^* E_G$.
This splitting $\beta\circ\beta'$ defines the tautological connection
on $\delta^* E_G$.

The above tautological connection on the
principal $G$--bundle $\delta^*E_G$ will be denoted by ${\mathcal D}_0$.

In Lemma \ref{lem2} we noted that the connections on $E_G$
are in bijective correspondence with the smooth sections
of ${\mathcal C}(E_G)$. Take any smooth section
\begin{equation}\label{chi}
\sigma\, :\, M\, \longrightarrow\, {\mathcal C}(E_G)
\end{equation}
of the fiber bundle ${\mathcal C}(E_G)
\, \longrightarrow\,M$. Let $D(\sigma)$
be the corresponding connection on the principal
$G$--bundle $E_G$.
We note that $\sigma^*\delta^*E_G\, =\, E_G$ because
$\delta\circ\sigma\, =\, \text{Id}_M$.

The following lemma is a consequence of
the construction of the tautological connection ${\mathcal D}_0$.

\begin{lemma}\label{lem0}
The connection $D(\sigma)$ on $E_G$ coincides with the
pulled back connection $\sigma^*{\mathcal D}_0$ on the
principal $G$--bundle $\sigma^*\delta^*E_G\, =\, E_G$.
\end{lemma}

\subsection{Construction of universal connection}

All infinite dimensional manifolds will be
modeled on the direct limit ${\mathbb R}^\infty$
of the sequence of vector spaces $\{{\mathbb R}^n\}_{n>0}$
with natural inclusions ${\mathbb R}^i\, \hookrightarrow\,
{\mathbb R}^{i+1}$.

Let
\begin{equation}\label{x0}
p_0\, :\, E_G\, \longrightarrow\, B_G
\end{equation}
be a universal principal $G$--bundle in the $C^\infty$ category;
see \cite{Mi} for the construction of a universal principal $G$--bundle.
So, $B_G$ is a $C^\infty$ manifold, the
projection $p_0$ is smooth, and
$E_G$ is contractible. Define
$$
{\mathcal B}_G\, :=\, B_G\times {\mathbb R}^\infty\, .
$$
Define
$$
{\mathcal E}_G\, :=\, p^*_{B_G} E_G\, =\, E_G\times{\mathbb R}^\infty\, ,
$$
where $p_{B_G}\, :\, B_G\times {\mathbb R}^\infty\,\longrightarrow
\, B_G$ is the natural projection.

See \cite[p.~374, Lemma~4.1]{BHS} for a proof of the following:

\begin{lemma}\label{lem3}
the principal $G$--bundle
$$
p\, :=\, p_0\times {\rm Id}_{{\mathbb R}^\infty}
\, :\, {\mathcal E}_G\,\longrightarrow\, {\mathcal B}_G
$$
is universal.
\end{lemma}

Set the principal $G$--bundle $E_G\, \longrightarrow\, M$ in Section
\ref{sec3} to be ${\mathcal E}_G\,\longrightarrow\, {\mathcal B}_G$.
Construct ${\mathcal C}({\mathcal E}_G)$ as
in \eqref{e11}. Let
\begin{equation}\label{x2}
\delta\, :\, {\mathcal C}({\mathcal E}_G)\, \longrightarrow\,
{\mathcal B}_G
\end{equation}
be the natural projection (see Lemma \ref{lem2}). Let ${\mathcal D}_0$
be the tautological connection on $\delta^*{\mathcal E}_G$ constructed
in Proposition \ref{prop2}.

The following theorem is proved in \cite[p.~375, Lemma~4.2]{BHS}.

\begin{theorem}\label{thm1}
The connection ${\mathcal D}_0$ on the principal
$G$--bundle $\delta^*{\mathcal E}_G$ is universal.
\end{theorem}

In Theorem \ref{thm1}, we took a special type of universal
$G$--bundle, namely we took the Cartesian product of a
universal $G$--bundle with ${\mathbb R}^\infty$. It should be
mentioned that Theorem \ref{thm1} is not valid if we do not
take this Cartesian product. For example, take $G$ to be the
additive group ${\mathbb R}^n$. Since ${\mathbb R}^n$ is contractible,
the projection
${\mathbb R}^n\, \longrightarrow\, \{\text{point}\}$ is a universal
${\mathbb R}^n$--bundle. Note that ${\mathcal C}({\mathbb R}^n)$
is a point. But the trivial principal ${\mathbb R}^n$
bundle on any manifold $X$ of dimension at least two admits connections
with nonzero curvature.

\section{Holomorphic connections}

Assume that $M$ is a complex manifold and $G$ is a complex Lie group. A {\it holomorphic}
principal $G$--bundle on $M$ is a triple $(E_G,\,  p,\, \psi)$ as in \eqref{tri} such that
$E_G$ is a complex manifold, and both the maps $p$ and $\psi$ are holomorphic.

Let $(E_G,\,  p,\, \psi)$ be a holomorphic principal $G$--bundle on $M$. Consider the
holomorphic tangent bundle $T^{1,0}E_G$, which is a holomorphic vector bundle on $E_G$.
The real tangent bundle $TE_G$ gets identified with $T^{1,0}E_G$ in the obvious way. More
precisely, the isomorphism $T^{1,0}E_G\, \longrightarrow\, TE_G$ sends a tangent vector
to its real part. Using this identification between $T^{1,0}E_G$ and $TE_G$, the complex
structure on the total space of $T^{1,0}E_G$ produces a complex structure on the total
space of $TE_G$. This complex structure on $TE_G$ produces a complex structure on the
quotient $\text{At}(E_G)$ in \eqref{e6}, because the action of $G$ on $TE_G$ is holomorphic.

The differential $dp$ in \eqref{dp} is holomorphic, which makes the projection
$\eta$ in \eqref{s.h.} holomorphic. The exact sequence in \eqref{e7} becomes an exact sequence of
holomorphic vector bundles. The holomorphic structure on $E_G$ produces a holomorphic structure
on any fiber bundle associated to $E_G$ for a holomorphic action of $G$. In particular, the
adjoint vector bundle $\text{ad}(E_G)$ has a holomorphic structure, because the
adjoint action of $G$ on $\mathfrak g$ is holomorphic. The homomorphism $\iota_0$
in \eqref{e7} is holomorphic with respect to this holomorphic structure on $\text{ad}(E_G)$.

A connection $$D\, :\, TM\,\longrightarrow\,{\rm At}(E_G)$$ on $E_G$ as in \eqref{D}
is called \textit{holomorphic} if the homomorphism $D$ is holomorphic.

\subsection{Holomorphic connection on principal bundles over a compact Riemann surface}

Now take $M$ to be a compact connected Riemann surface. It
is natural to ask the question when a holomorphic vector bundle on $M$
admits a holomorphic connection. Note that any holomorphic connection on
a Riemann surface is automatically flat because there are no nonzero
$(2, \,0)$ forms on a Riemann surface. A well-known theorem of
Atiyah and Weil says that a holomorphic vector
bundle $E$ over $M$ admits a holomorphic connection if and only
if each direct summand of $E$ is of degree zero (see \cite{At},
\cite{We}). We will describe a generalization of it to principal bundles.

Let $G$ be a complex connected reductive affine algebraic group. A parabolic subgroup of $G$ is
a Zariski closed connected subgroup $P\, \subset\, G$ such that the quotient
$G/P$ is compact. A Levi subgroup of of $G$ is a Zariski closed connected subgroup
$$
L\, \subset\, G
$$
such that there is a parabolic subgroup $P\, \subset\, G$ containing $L$ that satisfies
the following condition: $L$ contains a maximal torus of $P$, and moreover $L$ is a maximal
reductive subgroup of $P$. Given a holomorphic principal $G$--bundle $E_G$ on $M$
and a complex Lie subgroup $H\, \subset\, G$, a holomorphic reduction of $E_G$ to $H$
is given by a holomorphic section of the holomorphic fiber bundle $E_G/H$ over $M$. Let
$$q_H\, :\, E_G\,\longrightarrow\, E_G/H$$ be the quotient map. If $\nu\, :\, M\,
\longrightarrow\, E_G/H$ is a holomorphic section of the fiber bundle $E_G/H$, then
note that $q^{-1}_H(\nu(M))\, \subset\, E_G$ is a holomorphic principal $H$--bundle
on $M$. If $E_H$ is a holomorphic principal $H$--bundle on $M$, and $\chi$ is a holomorphic
character of $H$, then the associated holomorphic line bundle $E_H(\lambda)\, =\,
(E_H\times{\mathbb C})/H$ is the quotient of $E_H\times {\mathbb C}$, where $(z_1,\, c_1),\,
(z_2,\, c_2)\, \in\, E_H\times \mathbb C$ are identified if there is an element $g\,\in\, H$
such that
\begin{itemize}
\item $z_2\,=\, z_1g$, and

\item $c_2\,=\, \frac{c_1}{\lambda(g)}$.
\end{itemize}

The following theorem is proved in \cite{AB} (see
\cite[Theorem 4.1]{AB}).

\begin{theorem}\label{thmaz}
A holomorphic $G$--bundle $E_G$ over $M$ admits a holomorphic connection
if and only if for every triple of the form $(H,\, E_H,\, \lambda)$, where
\begin{enumerate}
\item{} $H$ is a Levi subgroup of $G$,
\item{} $E_H\subset E_G$ is a holomorphic reduction of structure
group to $H$, and
\item{} $\lambda$ is a holomorphic character of $H$,
\end{enumerate}
the associated line bundle $E_H(\lambda)\, =\,
(E_H\times{\mathbb C})/H$ over $M$ is of degree zero.
\end{theorem}

Note that setting $G\, =\, \text{GL}(n, {\mathbb C})$ in Theorem \ref{thmaz}
the above mentioned criterion of Atiyah and Weil is recovered.

We will describe a sketch of the proof of Theorem \ref{thmaz}.

Let $E_G$ be a
holomorphic $G$--bundle over $M$ equipped with a holomorphic
connection $\nabla$. Take any triple $(H,\, E_H,\, \lambda)$
as in Theorem \ref{thmaz}. We will first show that the connection $\nabla$ produces
a holomorphic connection on the principal $H$--bundle $E_H$.

Let $\mathfrak g$ and $\mathfrak h$ denote the Lie algebras of $G$ and $H$ respectively.
The group $H$ has adjoint actions on both $\mathfrak h$ and $\mathfrak g$.
To construct the connection on $E_H$, fix a splitting of the
injective homomorphism of $H$--modules
$$
0\, \longrightarrow\, {\mathfrak h}\, \longrightarrow\,
{\mathfrak g}\, .
$$
Since a holomorphic connection on $E_G$ is a given by a holomorphic splitting
of the Atiyah exact sequence for $E_G$, a holomorphic connection $\nabla$ on $E_G$ produces
a $\mathfrak g$--valued holomorphic $1$--form $\omega$ on $E_G$ satisfying the
following two conditions:
\begin{itemize}
\item $\omega$ is $G$--equivariant ($G$ acts on $\mathfrak g$ by inner automorphism), and

\item the restriction of $\omega$ to any fiber of $E_G$ is the Maurer--Cartan form on the
fiber.
\end{itemize}
Using the chosen splitting homomorphism
$${\mathfrak g}\, \longrightarrow\, {\mathfrak h}
\, \longrightarrow\, 0\, ,$$ the
connection form $\omega$ on $E_G$ defines a $\mathfrak h$--valued holomorphic
one--form $\omega'$ on $E_G$. The restriction of $\omega'$ to the complex
submanifold $E_H\, \subset\, E_G$ satisfies the two conditions needed for a holomorphic
$\mathfrak h$--valued $1$--form on $E_H$ to define a holomorphic connection on $E_H$.

Therefore, $E_H$ admits a holomorphic connection. A
holomorphic connection on $E_H$ induces a holomorphic connection
on the associated line bundle $E_H(\lambda)$. Any line
bundle admitting a holomorphic connection must be of degree
zero \cite{At}. Therefore, if $E_G$ admits a holomorphic connection then we know that
the degree of $E_H(\lambda)$ is zero.

To prove the converse, let $E_G$ be a
holomorphic $G$--bundle over $M$ such that
$$
\text{degree}(E_H(\lambda))\, =\,0
$$
for all triples
$(H,\, E_H,\, \lambda)$ of the above type. We need to show that
the Atiyah exact sequence for $E_G$ in \eqref{e7} splits holomorphically.

As the first step, in \cite{AB} the following is proved: it is enough to prove
that the Atiyah exact sequence for $E_G$ splits holomorphically under the
assumption that $E_G$ does not admit any
holomorphic reduction of structure group to any proper
Levi subgroup of $G$. Therefore, we assume that $E_G$ does not admit any
holomorphic reduction of structure group to any proper
Levi subgroup of $G$.

Let $\Omega^1_M$ denote the holomorphic cotangent bundle of $M$.
The obstruction for splitting of the Atiyah exact
sequence for $E_G$ is an element
$$
\tau(E_G)\, \in\, H^1(M,\, \Omega^1_M\otimes \text{ad}(E_G))\, .
$$
By Serre duality,
$$
H^1(M,\, \Omega^1_M\otimes \text{ad}(E_G))\,=\,
H^0(M,\, \text{ad}(E_G))^*\, .
$$
So we have
\begin{equation}\label{tau}
\tau(E_G)\, \in\, H^1(M,\, \text{ad}(E_G))^*\, .
\end{equation}

Any homomorphic section $f$ of $\text{ad}(E_G)$ has a Jordan
decomposition
$$
f\, =\, f_s+ f_n\, ,
$$
where $f_s$ is pointwise
semisimple and $f_n$ is pointwise nilpotent. From the
assumption that $E_G$ does not admit any
holomorphic reduction of structure group to any proper
Levi subgroup of $G$ it follows that the semisimple section $f_s$ is given by
some element of the center of $\mathfrak g$. Using this, from
the assumption on $E_G$ it can be deduced that
$$
\tau(E_G) (f_s) \, =\, 0\, ,
$$
where $\tau(E_G)$ is the element in \eqref{tau}.

The nilpotent section $f_n$ of $\text{ad}(E_G)$ gives a holomorphic reduction
of structure group $E_P\, \subset\, E_G$
of $E_G$ to a proper parabolic subgroup $P$ of $G$. This
reduction $E_P$ has the property that $f_n$ lies in the image
$$
H^0(M,\, \text{ad}(E_P))\, \hookrightarrow\,
H^0(M,\, \text{ad}(E_G))\, ,
$$
where $\text{ad}(E_P)$ is the adjoint bundle of $E_P$.
Using this reduction it can be shown that $\tau(E_G)(f_n)\, =\, 0$.

Hence $\tau(E_G)(f)\, =\, 0$ for all $f$, which implies
that $\tau(E_G)\, =\, 0$. Therefore, the Atiyah exact sequence for
$E_G$ splits holomorphically, implying that $E_G$ admits a holomorphic
connection.

\section{Real Higgs bundles}

As before, let $M$ be a compact connected Riemann surface. Let $$\sigma\,:\,
M\, \longrightarrow\, M$$ be an anti-holomorphic automorphism of order
two. Take a holomorphic vector bundle $E$ on $M$ of rank $r$. Let $\overline{E}$
denote the $C^\infty$ $\mathbb C$--vector bundle on $M$ of rank $r$ whose underlying
$C^\infty$ $\mathbb R$--vector bundle is the $\mathbb R$--vector bundle underlying $E$,
while the multiplication by $\sqrt{-1}$ on the fibers of $\overline{E}$ coincides with the
multiplication by $-\sqrt{-1}$ on the fibers of $E$. We note that the pullback
$\sigma^*\overline{E}$ has a natural structure of a holomorphic vector bundle. Indeed,
a $C^\infty$ section $s$ of $\sigma^*\overline{E}$ defined over an open subset
$U\, \subset\, M$ is holomorphic if the section $\sigma^*s$ of $E$ over $\sigma(U)$
is holomorphic; this condition uniquely defines the holomorphic structure on
$\sigma^*\overline{E}$. We use the terminology ``$\mathbb R$--vector bundles'' because
the terminology ``real vector bundles'' will be used for something else.

If $\alpha\, :\, A\, \longrightarrow\, B$ is a $C^\infty$ homomorphism of
holomorphic vector bundles on $M$, then $\overline{\alpha}$ will denote the homomorphism
$\overline{A}\, \longrightarrow\, \overline{B}$ defined by $\alpha$ using the
identifications of $A$ and $B$ with $\overline{A}$ and $\overline{B}$ respectively. 
A \textit{real structure} on $E$ is a holomorphic isomorphism of vector bundles
$$
\phi\, :\, E\, \longrightarrow\,\sigma^*\overline{E}
$$
over the identity map of $M$ such that the composition
\begin{equation}\label{rq}
E\, \stackrel{\phi}{\longrightarrow}\,\sigma^*\overline{E}\,
\stackrel{\sigma^*\overline{\phi}}{\longrightarrow}\,
\sigma^*\overline{\sigma^*\overline{E}}\,=\, E
\end{equation}
is the identity map of $E$.

A \textit{quaternionic structure} on $E$ is a holomorphic isomorphism of vector bundles
$$
\phi\, :\, E\, \longrightarrow\,\sigma^*\overline{E}
$$
over the identity map of $M$ such that the composition $E\, \longrightarrow\, E$
in \eqref{rq} is $-\text{Id}_E$.

A \textit{real vector bundle} on $(M,\, \sigma)$ is a pair of the form $(E,\, \phi)$, where
$E$ is a holomorphic vector bundle on $M$ and $\phi$ is a real structure on $E$.

A \textit{quaternionic vector bundle} on $(M,\, \sigma)$ is a pair of the form $(E,\, \phi)$, where
$E$ is a holomorphic vector bundle on $M$ and $\phi$ is a quaternionic structure on $E$.

Consider the differential $d\sigma\, :\, T^{\mathbb R}M \,\longrightarrow\,
\sigma^*T^{\mathbb R}M$ of the automorphism $\sigma$. Since $\sigma$ is anti-holomorphic,
it produces an isomorphism
$$
\sigma''\, :\, T^{1,0}M \,\longrightarrow\,
\sigma^*T^{0,1}M\,=\, \sigma^*\overline{T^{1,0}M}\, .
$$ 
It is easy to check that $\sigma''$ is holomorphic and it is a real structure on the
holomorphic tangent bundle $T^{1,0}M$. Let
\begin{equation}\label{el1}
\sigma'\, :\, K_M\, :=\, (T^{1,0}M)^* \,\longrightarrow\, \sigma^*\overline{K_M}
\end{equation}
be the real structure on the holomorphic cotangent bundle $K_M$ obtained from $\sigma''$.

We recall that a Higgs field on $E$ is a holomorphic section of
$\text{Hom}(E,\, E\otimes K_M)\,=\, \text{End}(E)\otimes K_M$ \cite{Hi},
\cite{Si}. A Higgs field $\theta$ on a real or quaternionic
vector bundle $(E,\, \phi)$ is called \textit{real} if the following diagram
is commutative:
$$
\begin{matrix}
E & \stackrel{\theta}{\longrightarrow} & E\otimes K_M\\
~\Big\downarrow\phi && \,~\,~\,~\,~\, ~\,\Big\downarrow\phi\otimes \sigma'\\
\sigma^*\overline{E}& \stackrel{\sigma^*\overline{\theta}}{\longrightarrow} &
\sigma^*\overline{E\otimes K_M}\,=\, \sigma^*\overline{E}\otimes \sigma^*\overline{K_M}
\end{matrix}
$$
where $\sigma'$ is the isomorphism in \eqref{el1}. A \textit{real} (respectively, quaternionic)
Higgs bundle on $(M,\, \sigma)$ is a triple of the form $((E,\, \phi),\, \theta)$, where
$(E,\, \phi)$ is a real (respectively, quaternionic) vector bundle on $(M,\, \sigma)$ and
$\theta$ is a real Higgs field on $(E,\, \phi)$.

We recall that the \textit{slope} of a holomorphic vector bundle $W$ on $M$ is the rational
number $\text{degree}(W)/\text{rank}(W)\,:=\, \mu(W)$. 
A real or quaternionic Higgs bundle $((E,\, \phi),\, \theta)$ on $(M,\, \sigma)$ is called 
\textit{semistable} (respectively, \textit{stable}) if for all nonzero holomorphic subbundle
$F\, \subsetneq \, E$ with
\begin{enumerate}
\item $\phi(F)\, \subset\, \sigma^*\overline{F}\, \subset\,\sigma^*\overline{E}$, and

\item $\theta(F)\, \subset\, F\otimes K_M$,
\end{enumerate}
we have
$\mu(F)\, \leq\, \mu(E)$ (respectively, $\mu(F)\, <\, \mu(E)$). A semistable
real (respectively, quaternionic) Higgs bundle is called \textit{polystable} if it is a
direct sum of stable real (respectively, quaternionic) Higgs bundles.

It is known that a real Higgs bundle $((E,\, \phi),\, \theta)$ is semistable (respectively, 
polystable) if and only if the Higgs bundle $(E,\, \theta)$ is
semistable (respectively, polystable) \cite[p.~2555, Lemma 5.3]{BGH}. Similarly,
a quaternionic Higgs bundle $((E,\, \phi),\, \theta)$ is semistable (respectively,     
polystable) if and only if the Higgs bundle $(E,\, \theta)$ is
semistable (respectively, polystable).

A polystable Higgs vector bundle $(E,\, \theta)$ of degree zero on $M$ admits a 
harmonic metric $h$ that satisfies the Yang--Mills--Higgs equation \cite{Si}, \cite{Do},
\cite{Hi}. If $((E,\, 
\phi),\, \theta)$ is real or quaternionic polystable of degree zero, then $E$ admits a harmonic
metric  $h$ because $(E,\, \theta)$ is polystable of degree zero. The harmonic metric $h$ on 
$E$ can be so chosen that the isomorphism $\phi$ is an isometry (note that $h$ induces 
a Hermitian structure on $\overline{E}$) \cite[p.~2557, Proposition 5.5]{BGH}.


\end{document}